\documentclass[12pt]{amsart}
\usepackage{amsmath,amsthm,amsfonts,amssymb}
\date{}
\newtheorem{theorem}{Theorem}
\newtheorem{lemma}{Lemma}
\newtheorem{corollary}{Corollary}
\newtheorem{proposition}{Proposition}

\newtheorem{problem}{Problem}

\newtheorem{example}{Example}
\theoremstyle{definition}\newtheorem{remark}{Remark}

\newcommand{\K}{\mathcal{K}}
\newcommand{\G}{\mathcal{G}}
\newcommand{\IN}{{\mathbb{N}}}
\newcommand{\IQ}{{\mathbb{Q}}}
\newcommand{\IC}{{\mathbb{C}}}
\newcommand{\IZ}{{\mathbb{Z}}}
\newcommand{\IT}{{\mathbb{T}}}
\newcommand{\IR}{{\mathbb{R}}}
\newcommand{\TopGr}{{\mathrm{TopGr}}}
\newcommand{\TBTG}{{\mathrm{TBG}}}
\newcommand{\TPTG}{{\mathrm{TPTG}}}
\newcommand{\FCTBTG}{{\mathrm{FCTBG}}}

\newcommand{\Arg}{\mathrm{Arg}}
\newcommand{\w}{\omega}

\newcommand{\bs}{\backslash}
\newcommand{\ol}{\overline}
\newcommand{\0}{\emptyset}
\newcommand{\inte}{\mathop{\rm int}\nolimits}
\renewcommand{\phi}{\varphi}

\begin{document}
\title[On subgroups of saturated paratopological groups]{On subgroups of saturated or totally bounded
paratopological groups}
\author{Taras Banakh}
\email{tbanakh@franko.lviv.ua}
\address{
Instytut Matematyki, Akademia \'Swi\c etokrzyska in Kielce,
\'Swi\c etokrzyska 15, Kielce, 25406, Poland\newline and\newline
Department of Mathematics, Ivan Franko Lviv National University,
Universytetska, 1 Lviv 79000, Ukraine.}
\author{Sasha Ravsky}
\email{oravsky@mail.ru}
\address{Department of Mathematics, Ivan Franko Lviv National University,
Universytetska, 1 Lviv 79000, Ukraine.}
\thanks{The first author was supported in part by the Slovenian-Ukrainian
research grant SLO-UKR 02-03/04.}


\hskip20pt\vbox{\hsize11truecm
\begin{abstract}
A paratopological group $G$ is {\em
saturated} if the inverse $U^{-1}$ of each non-empty set $U\subset
G$ has non-empty interior. It is shown that a [first-countable]
paratopological group $H$ is a closed subgroup of a saturated
(totally bounded) [abelian] paratopological group if and only if
$H$ admits a continuous bijective homomorphism onto a (totally
bounded) [abelian] topological group $G$ [such that for each
neighborhood $U\subset H$ of the unit $e$ there is a closed subset
$F\subset G$ with $e\in h^{-1}(F)\subset U$]. As an application we
construct a paratopological group whose character exceeds its
$\pi$-weight as well as the character of its group reflexion. Also
we present several examples of (para)topological groups which are
subgroups of totally bounded paratopological groups but fail to be
subgroups of {\em regular\/} totally bounded paratopological
groups.
\end{abstract}
}
\maketitle \baselineskip15pt

In this paper we continue investigations of paratopological
groups, started by the authors in \cite{Ra1}, \cite{Ra2}, \cite{Ra4},
\cite{BR1}--\cite{BR3}. By a {\em paratopological group} we
understand a pair $(G,\tau)$ consisting of a group $G$ and a
topology $\tau$ on $G$ making the group operation $\cdot:G\times
G\to G$ of $G$ continuous (such a topology $\tau$ will be called a
{\it semigroup topology} on $G$). If, in addition, the operation
$(\cdot)^{-1}:G\to G$ of taking the inverse is continuous with
respect to the topology $\tau$, then $(G,\tau)$ is a {\em
topological group\/}. All topological spaces considered in this
paper are supposed to be Hausdorff if the opposite is not stated.

The absence of the continuity of the inverse in paratopological groups results
in appearing various pathologies impossible in the category of topological
groups, which makes the theory of paratopological groups quite interesting and
unpredictable. In \cite{Gu} I.~Guran has introduced a relatively narrow class
of so-called saturated paratopological groups which behave much like usual
topological groups. Following I.Guran we define a paratopological group $G$ to
be {\em saturated\/} if the inverse $U^{-1}$ of any nonempty open subset $U$ of
$G$ has non-empty interior in $G$. A standard example of a saturated
paratopological group with discontinuous inverse is the {\em Sorgenfrey line},
that is the real line endowed with the Sorgenfrey topology generated by the
base consisting of half-intervals $[a,b)$, $a<b$. Important examples of
saturated paratopological groups are {\em totally bounded} groups, that is
paratopological groups $G$ such that for any neighborhood $U\subset G$ of the
origin in $G$ there is a finite subset $F\subset G$ with $G=FU=UF$, see \cite
[Proposition 3.1]{Ra3}.

Observing that each subgroup of the Sorgenfrey line is saturated, I.Guran asked
if the same is true for any saturated paratopological group. This question can
be also posed in another way: which paratopological groups embed into saturated
ones? A similar question concerning embedding into totally bounded semi- or
paratopological groups appeared in \cite{Ar}.

In this paper we shall show that the necessary and sufficient
condition for a paratopological group $G$ to embed into a
saturated (totally bounded) paratopological group is the existence
of a bijective continuous group homomorphism $h:G\to H$ onto a
topological (totally bounded) group $H$. The latter property of
$G$ will be referred to as the $\flat$-separateness (resp. Bohr
separateness).

$\flat$-Separated paratopological groups can be equivalently
defined with help of the group reflexion of a paratopological
group. Given a paratopological group $G$ let $\tau^\flat$ be the
strongest group topology on $G$, weaker than the topology of $G$.
The topological group $G^\flat=(G,\tau^\flat$), called {\em the
group reflexion\/} of $G$, has the following characteristic
property: the identity map $i:G\to G^\flat$ is continuous and for
every continuous group homomorphism $h:G\to H$ from $G$ into a
topological group $H$ the homomorphism $h\circ i^{-1}:G^\flat\to
H$ is continuous. According to \cite{BR1}, a neighborhood base of
the unit of the group reflexion $G^\flat$ of a {\em saturated}
(more generally, 2-oscillating) paratopological group $G$ consists
of the sets $UU^{-1}$ where $U$ runs over neighborhoods of the
unit in $G$. For instance, the group reflexion of the Sorgenfrey
line is the usual real line.

There is also a dual notion of a group co-reflexion. Given a
paratopological group $G$ let $\tau_\sharp$ be the weakest group
topology on $G$, stronger than the topology of $G$. The
topological group $G^\sharp=(G,\tau_\sharp)$ is called {\em the
group co-reflexion\/} of $G$. According to \cite{Ra1}, a
neighborhood base of the unit of the group co-reflexion $G^\sharp$
of a paratopological group $G$ consists of the sets $U\cap U^{-1}$
where $U$ runs over neighborhoods of the unit in $G$.  A
paratopological group is called {\em $\sharp$-discrete} provided
its group co-reflection is discrete. For instance, the Sorgenfrey
line is $\sharp$-discrete.

A subset $A$ of a paratopological group $G$ will be called {\em
$\flat$-closed\/} if $A$ is closed in the topology $\tau^\flat$. A
paratopological group $G$ is called {\em $\flat$-separated\/}
provided its group reflexion $G^\flat$ is Hausdorff; $G$ is called
{\em $\flat$-regular\/} if each neighborhood $U$ of the unit $e$
of $G$ contains a $\flat$-closed neighborhood of $e$. It is easy
to see that each $\flat$-regular paratopological group is regular
and  $\flat$-separated, see \cite{BR1}. For saturated groups the
converse is also true: each regular saturated paratopological
group is $\flat$-regular, see \cite[Theorem 3]{BR1}.

The notions of a $\flat$-separated (resp. $\flat$-regular)
paratopological group is a partial case of a more general notion
of a $\G$-separated (resp. $\G$-regular) paratopological group,
where $\G$ is a class of topological groups.  Namely, we call a
paratopological group $G$ to be
\begin{itemize}
\item {\em $\G$-separated} if $G$ admits a continuous bijective homomorphism
$h:G\to H$ onto a topological group $H\in\G$;
\item {\em $\G$-regular} if $G$ admits a regular continuous homomorphism
$h:G\to H$ onto a topological group $H\in\G$.
\end{itemize}

A continuous map $h:X\to Y$ between topological spaces is defined
to be {\em regular\/} if for each point $x\in X$ and a
neighborhood $U$ of $x$ in $X$ there is a closed subset $F\subset
Y$ such that $h^{-1}(F)$ is a closed neighborhood of $x$ with
$h^{-1}(F)\subset U$.

As we shall see later, any injective continuous map from a
$k_\omega$-space is regular. 
We remind that a topological space $X$ is defined to be a {\em
$k_\omega$-space} if there is a countable cover $\K$ of $X$ by
compact subsets of $X$, determining the topology of $X$ in the
sense that a subset $U$ of $X$ is open in $X$ if and only if the
intersection $U\cap K$ is open in $K$ for any compact set
$K\in\K$. According to \cite{FT}, each Hausdorff $k_\omega$-space
is normal. Under a {\em (para)topological $k_\omega$-group} we
shall understand a (para)topological group whose underlying
topological space is $k_\omega$. Many examples of
$k_\omega$-spaces appear in topological algebra as free objects in
various categories, see \cite{Cho}. In particular, the free
(abelian) topological group over a compact Hausdorff space is a
topological $k_\omega$-group, see \cite{Gra}.

\begin{proposition}\label{kos} Any injective continuous map $f:X\to Y$ from a
Hausdorff $k_\omega$-space $X$ into a Hausdorff space $Y$ is
regular.
\end{proposition}

\begin{proof} Fix any point $x\in X$
and an open neighborhood $U\subset X$ of $x$. Let $\{K_n\}$ be an
increasing sequence of compacta determining the topology of the
space $X$. Without loss of generality, we may assume that
$K_0=\{x\}$. Let $V_0=K_0$ and $W_0=Y$.

By induction, for every $n\ge 1$ we shall find an open
neighborhood $V_n$ of $x$ in $K_n$ and an open neighborhood $W_n$
of $f(\overline{V}_n)$ in $Y$ such that

1) $f^{-1}(\ol{W}_n)\cap K_n\subset U$;

2) $\ol{V}_n\subset U\cap \bigcap_{i\le n}f^{-1}(W_i)$;

3) $\ol{V}_n\subset V_{n+1}$.

Assume that for some $n$ the sets $V_i,W_i$, $i<n$, have been
constructed. It follows from (2) that $f(\ol{V}_{n-1})$ and
$f(K_{n}\setminus U)$ are disjoint compact sets in the Hausdorff
space $Y$. Consequently, the compact set $f(\ol{V}_{n-1})$ has an
open neighborhood $W_n\subset Y$ whose closure $\ol{W}_n$ misses
the compact set $f(K_n\setminus U)$. Such a choice of $W_n$ and
the condition (2) imply that  $\ol{V}_{n-1}\subset U\cap
\bigcap_{i\le n}f^{-1}(W_i)$. Using the normality of the compact
space $K_n$ find an open neighborhood $V_n$ of $\ol{V}_{n-1}$ in
$K_n$ such that $\ol{V}_n \subset U\cap \bigcap_{i\le
n}f^{-1}(W_i)$, which finishes the inductive step.

It is easy to see that $V=\bigcup_{n\in\w}V_n$ is an open
neighborhood of $x$ such that $f^{-1}(\ol{f(V)})\subset
f^{-1}(\bigcap_{i\in\w}\ol{W}_i)\subset U$, which just yields the
regularity of the map $f$.
\end{proof}

For paratopological $k_\omega$-groups this proposition yields the equivalence
between the $\G$-regularity and $\G$-separateness.

\begin{corollary}\label{kog} Let $\G$ be a class of topological groups.
A paratopological $k_\omega$-group is $\G$-separated if and only if it is
$\G$-regular.\qed
\end{corollary}

A class $\G$ of topological groups will be called
\begin{itemize}
\item {\em closed-hereditary} if each closed subgroup of a
group $G\in\G$ belongs to the class $\G$;
\item {\em $H$-stable}, where $H$ is a topological group, if $G\times
H\in\G$ for any topological group $G\in\G$.
\end{itemize}

For a topological space $X$ by $\chi(X)$ we denote its {\em
character}, equal to the smallest cardinal $\kappa$ such that each
point $x\in X$ has a neighborhood base of size $\le\kappa$.

 Now we are able to give a characterization of
subgroups of saturated paratopological groups (possessing certain
additional properties).

\begin{theorem}\label{t1} Suppose $T$ is a saturated $\sharp$-discrete
nondiscrete paratopological group and $\G$ is a closed-hereditary
$T^\flat$-stable class of topological groups. A paratopological
group $H$ is $\G$-separated if and only if $H$ is a $\flat$-closed
subgroup of a saturated paratopological group $G$ with
$G^\flat\in\G$, $\chi(G)=\max\{\chi(H),\chi(T)\}$, and
$|G/H|=|T|$.
\end{theorem}

A similar characterization holds for $\G$-regular paratopological
groups. We remind that a paratopological group $G$ is called a
{\em paratopological SIN-group} if any neighborhood $U$ of the
unit $e$ in $G$ contains a neighborhood $W\subset G$ of $e$ such
that $gWg^{-1}\subset U$ for all $g\in G$.  It is easy to check
that every paratopological SIN-group has a base at the unit
consisting of invariant open neighborhoods, see \cite[Ch.4]{Ra3}
(as expected, a subset $A$ of a group $G$ is called {\em
invariant} if $xAx^{-1}=A$ for all $x\in G$).

Finally we define a notion of a Sorgenfrey paratopological group
which crystallizes some important properties of the Sorgenfrey
topology on the real line. A paratopological group $G$ is defined
to be {\em Sorgenfrey\/} if $G$ is non-discrete, saturated and
contains a neighborhood $U$ of the unit $e$ such that for any
neighborhood $V\subset G$ of $e$ there is a neighborhood $W\subset
G$ of $e$ such that $x,y\in V$ for any elements $x,y\in U$ with
$xy\in W$. Observe that each Sorgenfrey paratopological group is
$\sharp$-discrete.

\begin{theorem}\label{t2} Suppose $T$ is a first countable saturated
regular Sorgenfrey paratopological  SIN-group and $\G$ is a
closed-hereditary $T^\flat$-stable class of first countable
topological SIN-groups. A paratopological group $H$ is
$\G$-regular if and only if $H$ is a $\flat$-closed subgroup of a
saturated $\flat$-regular paratopological group $G$ with
$G^\flat\in\G$, $\chi(G)=\chi(H)$, and $|G/H|=|T|$.
\end{theorem}

Theorem~\ref{t2} in combination with Corollary~\ref{kog} yield
\smallskip

\begin{corollary}~\label{cor2} Suppose $T$ is a first countable saturated regular Sorgenfrey
paratopological SIN-group and $\G$ is a closed-hereditary
$T^\flat$-stable class of first countable topological SIN-groups.
A paratopological $k_\omega$-group $H$ is $\G$-separated if and
only if $H$ is a $\flat$-closed subgroup of a saturated
$\flat$-regular paratopological group $G$ with $G^\flat\in\G$,
$\chi(G)=\chi(H)$, and $|G/H|=|T|$.\qed
\end{corollary}

As we said, any (regular) saturated paratopological group is
$\flat$-separated (and $\flat$-regular), see \cite{BR1}. Observe
that a paratopological group $G$ is $\flat$-separated if and only
if $G$ is $\TopGr$-separated where $\TopGr$ stands for the class
of all Hausdorff topological groups. These observations and
Theorem~\ref{t1} imply

\begin{corollary}\label{cor3} A paratopological group $H$ is $\flat$-separated
if and only if  $H$ is a ($\flat$-closed) subgroup of a saturated
paratopological group.\qed
\end{corollary}

Unfortunately we do not know the answer to the obvious
$\flat$-regular version of the above corollary.

\begin{problem} Is every $\flat$-regular paratopological group a
subgroup of a regular saturated paratopological group?
\end{problem}

For first-countable paratopological SIN-groups the answer to this
problem is affirmative.

\begin{corollary}\label{cor4} A first-countable paratopological SIN-group
$H$ is $\flat$-regular if and only if $H$ is a $\flat$-closed
subgroup of a regular first-countable saturated paratopological
SIN-group $G$ with $|G/H|=\aleph_0$.
\end{corollary}

\begin{proof} Taking into account that $H$ is a first-countable paratopological SIN-group
and applying \cite[Proposition 3]{BR1}, we conclude that $H^\flat$
is a first-countable topological SIN-group. Let $T$ be the
quotient group $\IQ/\IZ$ of the group of rational numbers, endowed
with the Sorgenfrey topology generated by the base consisting of
half-intervals. Observe that $T$ is a $\flat$-regular Sorgenfrey
abelian paratopological group and the class of first countable
topological SIN-groups is $T^\flat$ stable, closed-hereditary and
contains $H^\flat$. Now to finish the proof apply
Theorem~\ref{t2}.
\end{proof}

Next, we apply Theorems~\ref{t1} and \ref{t2} to the class $\TBTG$
(resp. $\FCTBTG$) of (first countable) totally bounded topological
groups. We remind that a paratopological group $G$ is {\em totally
bounded\/} if for any neighborhood $U$ of the unit in $G$ there is
a finite subset $F\subset G$ with $UF=FU=G$. It is known that each
totally bounded paratopological group is saturated and a saturated
paratopological group $G$ is totally bounded if and only if its
group reflexion $G^\flat$ is totally bounded, see \cite{Ra2},
\cite{BR1}. An example of a $\flat$-regular totally bounded
Sorgenfrey paratopological group is the circle
$\IT=\{z\in\IC:|z|=1\}$ endowed with the Sorgenfrey topology
generated by the base consisting of half-intervals $\{z\in\IT:a\le
\Arg(z)<b\}$ where $0\le a<b\le 2\pi$.

We define a paratopological group $G$ to be {\em Bohr separated}
(resp. {\em Bohr regular}, {\em fcBohr regular}) if it is
$\TBTG$-separated (resp. $\TBTG$-regular, $\FCTBTG$-regular). In
this terminology Theorem~\ref{t1} implies

\begin{corollary}\label{cor5} A paratopological group $H$ is Bohr separated if
and only if $H$ is a ($\flat$-closed) subgroup of a totally
bounded paratopological group.\qed
\end{corollary}

It is interesting to compare Corollary~\ref{cor5} with another
characterizing theorem supplying us with many pathological
examples of pseudocompact paratopological groups. We remind that a
topological space $X$ is {\em pseudocompact} if each locally
finite family of non-empty open subsets of $X$ is finite. It
should be mentioned that a Tychonoff space $X$ is pseudocompact if
and only if each continuous real-valued function on $X$ is
bounded.

\begin{theorem}\label{t3} An abelian paratopological group $H$ is Bohr separated
if and only if $H$ is a subgroup of a pseudocompact abelian
paratopological group $G$ with $\chi(G)=\chi(H)$.
\end{theorem}

The following characterization of fcBohr regular paratopological
groups follows from Theorem~\ref{t2} applied to the class
$\mathcal G=\FCTBTG$ and the quotient group $T=\IQ/\IZ$ endowed
with the standard Sorgenfrey topology.

\begin{corollary}\label{cor6} A paratopological group $H$ is
fcBohr regular if and only if $H$ is a $\flat$-closed subgroup of
a regular totally bounded paratopological group $G$ with
$\aleph_0=\chi(G^\flat)\le \chi(G)=\chi(H)$ and
$|G/H|\le\aleph_0$.
\end{corollary}

In some cases the fcBohr ($=\FCTBTG$) regularity is equivalent to
the Bohr ($=\TBTG$) regularity. We recall that a topological space
$X$ has {\em countable pseudocharacter} if each one point subset
of $X$ is a $G_\delta$-subset.

\begin{proposition}\label{p2} A Bohr regular paratopological group $H$ is
fcBohr regular provided one of the following conditions is
satisfied: \begin{enumerate}
\item $H$ is a $k_\omega$-space with countable pseudocharacter;
\item $H$ is first countable and Lindel\"of.
\end{enumerate}
\end{proposition}

\begin{proof} Using the Bohr regularity of $H$, find a
regular bijective continuous homomorphism $h:H\to K$ onto a
totally bounded topological group $K$. Denote by $e_H$ and $e_K$
the units of the groups $H,K$, respectively.

1. First assume that $H$ is a $k_\omega$-space with countable
pseudo-character. In this case the set $H\setminus\{e_H\}$ is
$\sigma$-compact as well as its image
$h(H\setminus\{e_H\})=K\setminus\{e_K\}$. It follows that the
totally bounded group $K$ has countable pseudo-character. Now it
is standard to find a bijective continuous homomorphism $i:K\to G$
of $K$ onto a first countable totally bounded topological group
$G$, see \cite[4.5]{Tk}. Since the composition $f\circ h:H\to G$
is bijective, the group $H$ is $\FCTBTG$-separated and by
Proposition~\ref{kos} is fcBohr regular.

2. Next assume that $H$ is first-countable and Lindel\"of. Fix a sequence
$(U_n)_{n\in\omega}$ of open subsets of $H$ forming a neighborhood base at
$e_H$. For every $n\in\omega$ fix a closed neighborhood $W_n\subset U_n$ whose
image $h(W_n)$ is closed in $K$. Let us call an open subset $U\subset K$ {\em
cylindrical} if $U=f^{-1}(V)$ for some continuous homomorphism $f:K\to G$ into
a first countable compact topological group $G$ and some open set $V\subset G$.
It follows from the total boundedness of $K$ that open cylindrical subsets form
a base of the topology of the group $K$, see \cite[3.4]{Tk}.

Using the Lindel\"of property of the set $f(H\setminus U_n)$, for every
$n\in\omega$ find a countable cover $\mathcal U_n$ of $f(H\setminus U_n)$ by
open cylindrical subsets such that $\cup\mathcal U_n\cap h(W_n)=\emptyset$.
Then $\mathcal U=\bigcup_{n\in\omega}\mathcal U_n$ is a countable collection of
open cylindrical subsets and we can produce a continuous homomorphism $f:K\to
G$ onto a first countable totally bounded topological group such that each set
$U\in\mathcal U$ is the preimage $U=f^{-1}(V)$ of some open set $V\subset K$.
To finish the proof it rests to observe that the composition $f\circ h:H\to G$
is a regular bijection of $H$ onto a first countable totally bounded
topological group.
\end{proof}

It can be shown that the character of any non-locally compact
paratopological $k_\omega$-group with countable pseudo-character
is equal to the small cardinal $\mathfrak d$, well-known in the
Modern Set Theory, see \cite{JW2}, \cite{Va}. By definition,
$\mathfrak d$ is equal to the cofinality of $\IN^\w$ endowed with
the natural partial order: $(x_i)_{i\in\w}\le (y_i)_{i\in\w}$ iff
$x_i\le y_i$ for all $i$. More precisely, $\mathfrak d$ is equal
to the smallest size of a subset $C\subset\IN^\w$ cofinal in the
sense that for any $x\in \IN^\w$ there is $y\in C$ with $y\ge x$.
It is easy to see that $\aleph_1\le\mathfrak d\le \mathfrak c$.
The Martin Axiom implies $\mathfrak d=\mathfrak c$. On the other
hand, there are models of ZFC with $\mathfrak d<\mathfrak c$, see
\cite{Va}.


We shall apply Corollary~\ref{cor6} to construct examples of
paratopological groups whose character exceed their $\pi$-weight
as well as the character of their group reflexions. We recall that
{\em the $\pi$-weight} $\pi w(X)$ of a topological space $X$ is
the smallest size of a {\em $\pi$-base}, i.e., a collection
$\mathcal B$ of open nonempty subsets of $X$ such that each
nonempty open subset $U$ of $X$ contains an element of the family
$\mathcal B$. According to \cite[4.3]{Tk} the $\pi$-weight of each
topological group coincides with its weight.

\begin{corollary}\label{cor7} For any uncountable cardinal $\kappa\le\mathfrak d$
there is a $\flat$-regular totally bounded countable abelian
paratopological group $G$ with $\aleph_0=\chi(G^\flat)=\pi
w(G)<\chi(G)=\kappa$.
\end{corollary}

\begin{proof} Take any non-metrizable countable abelian
$\FCTBTG$-separated topological $k_\omega$-group $(H,\tau)$ (for
example, let $H$ be a free abelian group over a convergent
sequence, see \cite{FT}). The group $H$, being
$\FCTBTG$-separated, admits a bijective continuous homomorphism
$h:H\to K$ onto a first countable totally bounded abelian
topological group $K$. By Proposition~\ref{kos} this homomorphism
is regular.

According to \cite[22)]{FT} or \cite{Ban}, the space $H$ contains
a copy of the Fr\'echet-Urysohn fan $S_\omega$ and thus has
character $\chi(H)\ge\chi(S_\omega)\ge\mathfrak d$. Using this
fact and the countability of $H$, by transfinite induction (over
ordinals $<\kappa$) construct a group topology $\sigma\subset
\tau$ on $H$ such that $\chi(H,\sigma)=\kappa$ and the
homomorphism $h:(H,\sigma)\to K$ is regular. This means that the
group $(H,\sigma)$ is fcBohr regular. Now we can apply
Corollary~\ref{cor6} to embed the group $(H,\sigma)$ into a
totally bounded countable abelian paratopological group $G$ such
that $\aleph_0=\chi(G^\flat)<\chi(G)=\chi(H,\sigma)=\kappa$. Since
$G$ is saturated and abelian, $\flat$-open subsets of $G$ form a
$\pi$-base which implies $\pi w(G)=\w(G^\flat)=\aleph_0$.
\end{proof}

\begin{remark} It is interesting to compare Corollary~\ref{cor7} with a result of
\cite{BRZ} asserting that there exists a $\flat$-regular countable
paratopological group $G$ with
$\aleph_0=\chi(G)<\chi(G^\flat)=\mathfrak d$. Such a
paratopological group $G$ cannot be saturated since
$\chi(G^\flat)\le\chi(G)$ for any saturated (more generally, any
2-oscillating) paratopological group $G$, see \cite{BR1}.
\end{remark}

We finish our discussion with presenting examples of regular
(para)topological groups which embed into totally bounded
paratopological groups but fail to embed into {\em regular}
totally bounded paratopological groups. For that it suffices to
find a Bohr separated group which is not Bohr regular.

Let us remark that each locally convex linear topological space
(or more generally each locally quasi-convex abelian topological
group, see \cite{A} or \cite{Ba}) is Bohr regular. On the other
hand, there exist  (non-locally convex) linear metric spaces which
fail to be Bohr separated or Bohr regular. The simplest example
can be constructed as follows. Consider the linear space $C[0,1]$
of continuous real-valued functions on the closed interval $[0,1]$
and endow it with the invariant metrics
$d_{1/2}(f,g)=\int_0^1\sqrt{|f(t)-g(t)|}dt$,
$p(f,g)=\sum_{n\in\w}\min\{2^{-n},|f(t_n)-g(t_n)|\}$ and
$\rho(f,g)=d_{1/2}(f,g)+p(f,g)$ where $\{t_n:n\in\w\}$ is an
enumeration of rational numbers of $[0,1]$. It is well-known that
the linear metric space $(C[0,1],d_{1/2})$ admits no non-zero
linear continuous functional and fails to be Bohr separated. The
linear metric space $(C[0,1],\rho)$ is even more interesting. We
remind that an abelian group $G$ is {\em divisible} (resp. {\em
torsion-free}) if for any $a\in G$ and natural $n$ the equation
$x^n=a$ has a solution (resp. has at most one solution) $x\in G$.

\begin{proposition}\label{lin}
The linear metric space $L=(C[0,1],\rho)$ is Bohr separated but
fails to be Bohr regular. Moreover, $L$ is a $\flat$-closed
subgroup of a totally bounded abelian torsion-free divisible
group, but fails to be a subgroup of a {\em regular} totally
bounded paratopological group.
\end{proposition}

\begin{proof} The Bohr separatedness of $L$ follows from the continuity
of the maps $\chi_n:L\to\IR$, $\chi_n:f\mapsto f(t_n)$, for
$n\in\omega$. Let us show that the group $L$ fails to be Bohr
regular.

For this we first prove that each linear continuous functional
$\psi:(L,\rho)\to\IR$
is continuous with respect to the ``product'' metric $p$. Consider
the open convex subset $C=\psi^{-1}(-1,1)$ of $L$. By the
continuity of $\psi$, there are $n\ge 1$ and $\varepsilon>0$ such
that $x\in C$ for any $x\in L$ with $d_{1/2}(x,0)<\varepsilon$ and
$|x(t_i)|<\varepsilon$ for all $i\le n$. Let $L_0=\{x\in
L:x(t_i)=0$ for all $i\le n\}$ and observe that the convex set
$C\cap L_0$ contains the open $\varepsilon$-ball with respect to
the restriction of the metric $d_{1/2}$ on $L_0$. Now the standard
argument (see \cite[4.16.3]{Ed}) yields $C\cap L_0=L_0$ and
$L_0\subset\bigcap_{k\ge1}\frac1kC=\psi^{-1}(0)$. Hence the
functional $\psi$ factors through the quotient space $L/L_0$ and
is continuous with respect to the metric $p$ (this follows from
the continuity of the quotient homomorphism $L\to L/L_0$ with
respect to $p$).

If $\chi:L\to\IT$ is any character on $L$ (that is a continuous
group homomorphism into the circle $\IT=\IR/\IZ$), then it is easy
to find a continuous linear functional $\psi:L\to \IR$ such that
$\chi=\pi\circ \psi$, where $\pi:\IR\to\IT$ is the quotient
homomorphism. As we have already shown, the functional $\psi$ is
continuous with respect to the metric $p$ and so is the character
$\chi$.

Finally, we are able to prove that the group $L$ fails to be Bohr
regular. Assuming the converse we would find a continuous regular
homomorphism $h:L\to H$ onto a totally bounded abelian topological
group $H$. The group $H$, being abelian and totally bounded, is a
subgroup of the product $\IT^\kappa$ for some cardinal $\kappa$,
see \cite{Po}. Then the above discussion yields that $h$ is
continuous with respect to the metric $p$. In this situation the
regularity of $h$ implies the regularity of the identity map
$(L,\rho)\to(L,p)$. But this map certainly is not regular: for any
$2^{-n}$-ball $B=\{x\in L:\rho(x,0)<2^{-n}\}$ its closure in the
metric $p$ contains the linear subspace $\{x\in L:x(t_i)=0$ for
all $i\le n\}$ and thus lies in no ball. Therefore the group $L$
is Bohr separated but not Bohr regular.

Let $\mathcal G$ be the class of all totally bounded abelian
divisible torsion-free topological groups. The group $L$, being
Bohr separated, abelian, divisible, and torsion-free, is $\mathcal
G$-separated. Pick any irrational number $\alpha\in\IT=\IR/\IZ$
and consider the subgroup $T=\{q\alpha:q\in\IQ\}$ of the circle
$\IT$ endowed with the Sorgenfrey topology. It is clear that $T$
is a totally bounded $\sharp$-discrete paratopological group with
$T^\flat\in\mathcal G$. By Theorem~\ref{t1}, $L$ is a
$\flat$-closed subgroup of a saturated paratopological groups $G$
with $G^\flat\in\mathcal G$ which implies that $G$ is totally
bounded abelian, divisible and torsion-free.

On the other hand, $L$ admits no embedding into a regular totally
bounded paratopological group $G$. Indeed, assuming that $L\subset
G$ is such an embedding, apply Theorem 3 of \cite{BR1} to conclude
that the identity homomorphism $id:G\to G^\flat$ is regular and so
is its restriction $id|L$, which would imply the Bohr regularity
of $L$.
\end{proof}

There is also an alternative method of constructing Bohr separated
but not Bohr regular paratopological groups, based on the concept
of a Lawson paratopological group. Following \cite{BR1} we define
a paratopological group $G$ to be {\em Lawson} if it has a
neighborhood base at the unit consisting of subsemigroups of $G$.
According to \cite{BR1} there is a regular Lawson paratopological
group failing to be $\flat$-separated. On the other hand, there
are Lawson paratopological groups which are $\flat$-regular and
Bohr separated but are not topological groups, see Example 2
\cite{BR1} or Example~\ref{Lawex} below. We shall show that a
$\flat$-regular paratopological group $G$ is a topological group
provided its group reflexion $G^\flat$ is topologically periodic.
We remind that a paratopological group $G$ is {\em topologically
periodic} if for each $x\in G$ and a neighborhood $U\subset G$ of
the unit there is a number $n\ge 1$ such that $x^n\in U$, see
\cite{BG}. It is easy to show that each totally bounded
topological group is topologically periodic. For paratopological
groups it is not true: according to Theorem~\ref{t2} there is a
$\flat$-regular totally bounded paratopological group $G$ which
contains the discrete group $\IZ$ of integers and thus cannot be
topologically periodic. The class of topologically periodic
topological groups will be denoted by $\TPTG$.

\begin{proposition}\label{Lawreg} Each $\TPTG$-regular Lawson
paratopological group is a topological group.
\end{proposition}

\begin{proof} Let $(G,\tau)$ be a Lawson paratopological group and
$\sigma\subset\tau$ be a topology turning $G$ into a topologically
periodic  topological group such that $(G,\tau)$ has a base
$\mathcal B$ at the unit consisting of subsemigroups, closed in
the topology $\sigma$. We are going to show that an arbitrary
element $U\in\mathcal B$ is in fact a subgroup of $G$. For this
purpose suppose that there exists an element $x\in U^{-1}\bs U$.
Then $x^{-1}\in U$ and $x^m\in U$ for all $m<0$ because $U$ is a
subsemigroup of $G$. Since the set $U$ is closed in the topology
$\sigma$, there exists a neighborhood $V\in\sigma$ of unit such
that $xV\cap U=\0$. By the topological periodicity of
$(G,\sigma)$, there exists a number $n<-1$ with $x^n\subset V$.
Then $x^{n+1}\cap U=\0$ which is a contradiction.
\end{proof}

Since each totally bounded topological group is topologically
periodic this Proposition implies

\begin{corollary}\label{LawBohreg} Each Bohr regular Lawson paratopological group
is a topological group.\qed
\end{corollary}

On the other hand, {\em abelian} Lawson paratopological groups are
Bohr separated.

\begin{proposition}\label{LawBohrsep} Each abelian Lawson
paratopological group is Bohr
separated.
\end{proposition}

\begin{proof} Let $G$ be such the group. Then $G^\flat$ has a
neighborhood base $\mathcal B$ at the unit, consisting of
subgroups. For every group $H\in\mathcal B$ the group $G/H$, being
abelian and discrete, is Bohr separated \cite{Po}. Since the
family $\{G\to G/H:H\in\mathcal B\}$ of quotient maps separates
the points of the group $G$, the group $G$ is Bohr separated too.
\end{proof}

Corollary~\ref{LawBohreg} and Proposition~\ref{LawBohrsep} allow
us to construct simple examples of Bohr separated Lawson
paratopological groups which are not Bohr regular.

\begin{example}\label{Lawex} There is a countable $\flat$-regular
saturated Lawson paratopological abelian group $H$ which is Bohr
separated but not Bohr regular. The group $H$ has the following
properties:
\begin{enumerate}
\item $H$ is a $\flat$-closed subgroup of a countable first-countable abelian totally bounded
paratopological group;
\item $H$ is a $\flat$-closed subgroup of a first-countable abelian pseudocompact
paratopological group;
\item $H$ fails to be a subgroup of a regular totally bounded (or
pseudocompact) paratopological group.
\end{enumerate}
\end{example}

\begin{proof} Consider the direct sum
$\IZ^\w_0=\{(x_i)_{i\in\omega}\in \IZ^\w: x_i=0$ for all but
finitely many indices $i\}$ of countably many copies of the group
$\IZ$ of integers. Endow the group $\IZ^\w_0$ with a shift
invariant topology $\tau$ whose neighborhood base at the origin
consists of the sets $U_n=\{0\}\cup \bigcup_{m\ge n}W_{m}$ where
$W_{m}=\{(x_i)_{i\in\omega}\in \IZ^\w_0: x_i=0$ for all $i<m$ and
$x_m>0\}$ for $m\ge 0$. It is easy to see that $H=(\IZ^\w_0,\tau)$
is a $\flat$-regular countable first-countable saturated Lawson
paratopological group which is not a topological group.
By Proposition~\ref{LawBohrsep} and Corollary~\ref{LawBohreg} the
group $H$ is Bohr separated but not Bohr regular.

By Theorem~\ref{t1}, $H$ is a $\flat$-closed subgroup of a
first-countable totally bounded countable paratopological group
and by Theorem~\ref{t3}, $H$ is a $\flat$-closed subgroup of a
first-countable abelian pseudocompact paratopological group.

Assuming that $H$ is a subgroup of a regular totally bounded or
pseudocompact paratopological group $G$ and applying Theorem 3 of
\cite{BR1} and \cite{RR} we would get that both $G$ and $H$ are
Bohr regular which is impossible.
\end{proof}

In the proofs of our principal results we shall often exploit the
following characterization of semigroup topologies on groups from
\cite[1.1]{Ra1}.

\begin{lemma}\label{Pon} A family $\mathcal B$ of subsets containing
a unit $e$ of a group $G$ is a neighborhood base at $e$ of some
semigroup topology $\tau$ on $G$ if and only if $\mathcal B$
satisfies the following four Pontryagin conditions:

1. $(\forall U,V\in\mathcal B)(\exists W\in\mathcal B):W\subset
U\cap V$;

2. $(\forall U\in\mathcal B)(\exists V\in\mathcal B):V^2\subset
U$;

3. $(\forall U\in\mathcal B)(\forall x\in U)(\exists V\in\mathcal
B):xV\subset U$;

4. $(\forall U\in\mathcal B)(\forall x\in G)(\exists V\in\mathcal
B):x^{-1}Vx\subset U$.

\noindent The topology $\tau$ is Hausdorff if and only if

5. $\bigcap\{UU^{-1}:U\in\mathcal B\}=\{e\}$.
\end{lemma}

\section{Proof of Theorem~\ref{t1}}

The necessity is evident. We shall prove the sufficiency. Let
$(H,\tau)$ be a $\mathcal G$-separated paratopological group,
where $\mathcal G$ is $T^\flat$-stable class of topological
groups. Since the group $H$ is $\mathcal G$-separated, there
exists a group topology $\sigma$ on the group $H$ such that
$(H,\sigma)\in\mathcal G$. We shall define the topology on the
product $G=H\times T$ as follows. Let $\mathcal B_\tau$, $\mathcal
B_\sigma$ and $\mathcal B_T$ be open bases at the unit of the
groups $(H,\tau)$, $(H,\sigma)$ and $T$ respectively. For
arbitrary neighborhoods $U_\tau\in \mathcal B_\tau$, $U_\sigma\in
\mathcal B_\sigma$ and $U_T\in \mathcal B_T$ with $U_\tau\subset
U_\sigma$ put  $[U_\tau, U_\sigma,U_T]=U_\tau\times\{e_T\}\cup
U_\sigma\times (U_T\bs\{e_T\})$, where $e_H$ and $e_T$ are the
units of the groups $H$ and $T$ respectively. The family of all
such $[U_\tau, U_\sigma,U_T]$ will be denoted by $\mathcal B$. Now
we verify the Pontryagin conditions for the family $\mathcal B$.

The Condition 1 is trivial.

To check Condition 2 consider an arbitrary set $[U_\tau,
U_\sigma,U_T]\in\mathcal B$. There exist neighborhoods
$V_\tau\in\mathcal B_\tau$, $V_\sigma\in\mathcal B_\sigma$ such
that $V_\tau^2\subset U_\tau$, $V_\sigma^2\subset U_\sigma$ and
$V_\tau\subset V_\sigma$. Since the group $T^\#$ is discrete then
there is a neighborhood $V_T\subset\mathcal B_T$ such that
$(V_T\bs\{e_T\})^2\subset U_T\bs\{e_T\}$. Then $[V_\tau,
V_\sigma,V_T]^2\subset [U_\tau, U_\sigma,U_T]$.

To verify Condition 3 consider an arbitrary point $x\in
[U_\tau,U_\sigma,U_T]\in\mathcal B$. If $x=(x_H,e_T)$, where
$x_H\in U_\tau$ then there exist neighborhoods $V_\tau\in\mathcal
B_\tau$, $V_\sigma\in\mathcal B_\sigma$ such that $V_\tau\subset
V_\sigma$, $x_HV_\tau\subset U_\tau$ and $x_HV_\sigma\subset
U_\sigma$. Then $x[V_\tau,V_\sigma, U_T]\subset [U_\tau,U_\sigma,
U_T]$. If $x=(x_H,x_T)$, where $x_H\in U_\sigma$ and $x_T\in
U_T\bs\{e_T\}$ then there exist neighborhoods $V_\tau\in\mathcal
B_\tau$, $V_\sigma\in\mathcal B_\sigma$ and $V_T\in\mathcal B_T$
such that $V_\tau\subset V_\sigma$, $x_HV_\sigma\subset U_\sigma$
and $x_TV_T\subset U_T\bs\{e_T\}$. Then $x[V_\tau, V_\sigma,
V_T]\subset [U_\tau, U_\sigma, U_T]$.

Condition 4. Let $x=(x_H,x_T)\subset H\times T$ be an arbitrary
point. Then there are neighborhoods $V_\tau\in\mathcal B_\tau$,
$V_\sigma\in\mathcal B_\sigma$ and $V_T\in\mathcal B_T$ such that
$V_\tau\subset V_\sigma$, $x_H^{-1}V_\tau x_H\subset U_\tau$,
$x_H^{-1}V_\sigma x_H\subset U_\sigma$ and $x_T^{-1}V_Tx_T\subset
U_T$. Then $x^{-1}[V_\tau,V_\sigma, V_T]x\subset [U_\tau,U_\sigma,
U_T]$.

Hence the family $\mathcal B$ is a base of a semigroup topology on
the group $G$. Denote this semigroup topology by $\rho$. The
inclusion
$\bigcap\{[U_\tau,U_\sigma,U_T]\cdot[U_\tau,U_\sigma,U_T]^{-1}:U_\tau\in\mathcal
B_\tau, U_\sigma\in\mathcal B_\sigma, U_T\in\mathcal
B_T\}\subset\{U_\sigma U_\sigma^{-1}\times
U_TU_T^{-1}:U_\sigma\in\mathcal B_\sigma, U_T\in\mathcal
B_T\}=\{(e_H,e_T)\}$ implies that the topology $\rho$ is
Hausdorff. Since the groups $T$ and $(H,\sigma)$ are saturated and
the group $T$ is nondiscrete, the group $(G,\rho)$ is saturated
too. According to \cite[Proposition 3]{BR1} the base at the unit
of the topology $\rho^\flat$ consists of the sets $UU^{-1}$, where
$U\in\mathcal B$. Thus the topology $\rho^\flat$ coincides with
the product topology of the groups $(H,\sigma)\times T^\flat$ and
hence $(G,\rho^\flat)\in\mathcal G$ and $H$ is a $\flat$-closed
subgroup of the group $G$.

\section{Proof of Theorem~\ref{t2}}

The ``if'' part of Theorem~\ref{t2} is trivial. To prove the
``only if'' part, suppose that $T$ and $(H,\tau)$ are
paratopological groups with the units $e_T$ and $e_H$, satisfying
the hypothesis of Theorem~\ref{t2}.

Using the Sorgenfrey property of the group $T$, choose an open
invariant neighborhood $U_0$ of the unit $e_T$ such that for any
neighborhood $U\subset T$ of $e_T$ there is a neighborhood
$U'\subset T$ of $e_T$ such that $x,y\in U$ for any elements
$x,y\in U_0$ with $xy\in U'$. By induction we can build a sequence
$\{U_n:n\in\omega\}$ of invariant open neighborhoods of $e_T$
satisfying the following conditions:

(1) $\{U_n:n\in\omega\}$ is a neighborhood base at the unit $e_T$
of the group $T$;

(2) $U_{n+1}^2\subset U_n$ for every $n\in\omega$;

(3) for every $n\in\omega$ and any points $x,y\in U_0$ the
inclusion  $xy\in U_{n+1}$ implies $x,y\in U_n$;

(4) $\ol {U_{n}}^\flat\subsetneqq U_{n-1}$ for every $n\in\omega$,
where $\ol {U_{n}}^\flat$ denotes the closure of the set $U_{n}$
in the topology of $T^\flat$.

Remark that the condition (3) yields

(5) $(U_0\bs U_n)U_0\cap U_{n+1}=\emptyset$ and hence $U_0\bs
U_n\cap U_{n+1}U_0^{-1}=\emptyset$ for all $n$.

\noindent Since the group $T$ is saturated, we can apply
Proposition~3 of \cite{BR1} to conclude that the set
$U_{n+2}U_0^{-1}$ is a neighborhood of the unit in $T^\flat$. Then
the set $U_{n+2}U_{n+2}U_0^{-1}\subset U_{n+1}U_0^{-1}$ is a
neighborhood of $U_{n+2}$ in $T^\flat$. This observation together
with (5) yields

(6) $\overline{U_0\bs U_n}^\flat\cap U_{n+2}=\emptyset$ for all
$n$.

It follows from our assumptions on $(H,\tau)$ that there exists a
group topology $\sigma\subset\tau$ on $H$ such that the group
$(H,\sigma)$ belongs to the class $\mathcal G$ and $(H,\tau)$ has
a neighborhood base $\mathcal B_\tau$ at the unit $e_H$ consisting
of sets, closed in the topology $\sigma$. By induction we can
build a base $\{V_n:n\in\omega\}$ of open symmetric invariant
neighborhoods of $e_H$ in the topology $\sigma$ such that
$V_{n+1}^2\subset V_n$ for every $n\in\omega$.

Consider the product $H\times T$ and identify $H$ with the
subgroup $H\times\{e_T\}$ of $H\times T$. It rests to define a
topology on $H\times T$. At first we shall introduce an auxiliary
sequence $\{W_k\}$ of ``neighborhoods'' of $(e_H,e_T)$ satisfying
the Pontryagin Conditions 1,2, and 4. For every $k\in\omega$ let
\medskip

$(\star)$\quad $ W_n=\{(e_H,e_T)\}\cup\bigcup\limits_{i>2n}
V_{ni}\times (U_{i-1}\setminus U_i)$ \smallskip

\noindent and observe that $W_{n+1}\subset W_n$ for all $n$. Let
us verify the Pontryagin Conditions 1,2,4 for the sequence
$(W_n)$.

To verify Conditions 1 and 2 it suffices to show that
$W_{n}^2\subset W_{n-1}$ for all $n\ge 1$. Fix any elements
$(x,t),(x',t')\in W_n$. We have to show that $(xx',tt')\in
W_{n-1}$. Without loss of generality, we can assume that $t,t'\ne
e_T$. In this case we may find numbers $i,i'> 2n$ with $(x,t)\in
V_{ni}\times (U_{i-1}\bs U_i)$ and $(x',t')\in
V_{ni'}\times(U_{i'-1}\bs U_{i'})$. For $j=\min\{i,i'\}$ the
Conditions (2), (5) imply
\begin{multline*} (xx',tt')\in
V_{nj-1}\times (U_{j-2}\bs U_{j+1})\subset
\bigcup_{k=j-1}^{j+1}V_{(n-1)k}\times (U_{k-1}\bs U_k)\subset \\
\bigcup_{k>2(n-1)}V_{(n-1)k}\times (U_{k-1}\bs U_k)\subset
W_{n-1}.\end{multline*}

Taking into account that both the sequences $\{U_n\}$ and
$\{V_n\}$ consist of invariant neighborhoods, we conclude that the
sets $W_n$ are invariant as well. Hence the Condition 4 holds too.

Now, using the sequence $(W_n)$ we shall produce a sequence
$(O_n)$ satisfying all the Pontryagin Conditions 1--5. For every
$n\in\omega$ put $O_n=\bigcup_{i=n}^\infty W_nW_{n+1}\cdots W_i$.
Thus $W_n\supset O_{n+1}\supset W_{n+1}$ and $O_n\cap
H\times\{e_T\}=\{(e_H,e_T)\}$ for all $n$. It is easy to see that
the sequence $\{O_n\}$ consists of invariant sets and satisfies
Pontryagin conditions 1--4. Hence the family $\{O_n\}$ is a
neighborhood base at the unit of some (not necessarily Hausdorff)
topology $\tau'$ on $G=H\times T$ turning $G$ into a
paratopological SIN-group. Applying Proposition 1.3 from
\cite{Ra1} we conclude that the family $\mathcal
B_\rho=\{OU:O\in\mathcal B_{\tau'}, U\in\mathcal B_\tau \}$ is a
neighborhood base at the unit of some (not necessarily Hausdorff)
semigroup topology $\rho$ on $G$ (here we identify $H$ with the
subgroup $H\times\{e_T\}$ in $G$). Since the topology $\rho$ is
stronger than the product topology $\pi$ of the group
$(H,\sigma)\times T^\flat$, the topology $\rho$ is Hausdorff and
$H$ is a $\flat$-closed subgroup of the group $(G,\rho)$. It
follows from the construction of the topology $\rho$ that
$\rho|H=\tau$, $\chi(G,\rho)=\chi(H)$ and $|G/H|=|T|$.

At the end of the proof we show that the paratopological group
$(G,\rho)$ is saturated and $\flat$-regular. To show that the
group $(G,\rho)$ is saturated it suffices to find for every $n\ge
1$ nonempty open sets $V\subset (H,\sigma)$ and $U\subset T$ such
that $V\times U^{-1}\subset W_n$. Taking into account that the
group $T$ is saturated and the set $U_{3n-1}\bs
\overline{U_{3n}}^\flat$ is nonempty, find a nonempty open set
$U\subset T$ such that $U^{-1}\subset U_{3n-1}\bs
\overline{U_{3n}}^\flat$. Then $V^{-1}_{3n^2}\times U^{-1}\subset
V_{3n^2}\times(U_{3n-1}\setminus U_{3n})\subset W_n$. This implies
that the group $(G,\rho)$ is saturated and
$(G,\rho^\flat)=(H,\sigma)\times T^\flat\in\mathcal G$.

The $\flat$-regularity of the group $(G,\rho)$ will follow as soon
as we prove that $\overline {W_{n}V}^\pi \subset W_{n-1}V$ for
every $n\ge 2$ and $V\in\mathcal B_\tau$. Indeed, in this case, we
shall get $$\overline{O_{n+1}V}^\flat\subset
\overline{O_{n+1}V}^\pi\subset\overline{W_nV}^\pi\subset
W_{n-1}V\subset O_{n-1}V.$$

Fix any $x\in\overline{W_nV}^\pi$. If $x\in V\times\{e_T\}$, then
$x\in W_{n-1}V$. Next, assume that $x\notin H\times\{e_T\}$. The
property (4) of the sequence $(U_k)$ implies that the point $x$
has a $\pi$-neighborhood meeting only finitely many sets $H\times
U_i$, $i\in\omega$. This observation together with $x\in
\overline{W_nV}^\pi$ and ($\star$) imply that
$x\in\overline{V_{ni}V\times (U_{i-1}\setminus U_i)}^\flat$ for
some $i>2n$. The condition (6) implies that the following chain of
inclusions holds:
\begin{multline*}
x\in\overline{V_{ni}V\times (U_{i-1}\setminus U_i)}^\flat\subset
\overline{V_{ni}V}^\sigma\times (\overline{U_{i-1}\setminus
U_i)}^\flat\subset V_{ni}^2V\times (U_{i-2}\setminus
U_{i+2})\subset\\ \bigcup_{j=i-1}^{i+2}
V_{ni-1}V\times(U_{j-1}\setminus U_j)\subset \bigcup_{j>
2n-2}V_{(n-1)j}V\times (U_{j-1}\setminus U_j)\subset W_{n-1}V.
\end{multline*}

Finally, assume that $x\in H\setminus V=(H\setminus
V)\times\{e_T\}$. Since the set $V$ is $\flat$-closed in $H$,
there is $m\in\omega$ such that $V_m^{-1}V_mx\cap V=\emptyset$ and
thus $V_mx\cap V_iV=\emptyset$ for all $i\ge m$. The inclusion
$x\in \overline{W_nV}^\pi$ and ($\star$) imply $$(V_m\times
U_mU_m^{-1})x\cap (V_{ni}V\times (U_{i-1}\setminus
U_i))\ne\emptyset$$ for some $i>2n$. Then $V_mx\cap
V_{ni}V\ne\emptyset$ and $U_mU_m^{-1}\cap (U_{i-1}\setminus
U_i)\ne\emptyset$. In view of Property (5) of the sequence
$(U_k)$, the latter relation implies $m\le i$. On the other hand,
the former relation together with the choice of the number $m$
yields $ni<m\le i$ which is impossible. This contradiction
finishes the proof of the inclusion $\overline{W_nV}^{\pi}\subset
W_{n-1}V$.

\section{Proof of Theorem~\ref{t3}}

Given a topological space $(X,\tau)$ Stone \cite{Sto} and Katetov
\cite{Kat} considered the topology $\tau_r$ on $X$ generated by
the base consisting of all canonically open sets of the space
$(X,\tau)$. This topology is called the {\it regularization} of
the topology $\tau$. If $(X,\tau)$ is Hausdorff then $(X,\tau_r)$
is regular and if $(X,\tau)$ is a paratopological group then
$(X,\tau_r)$ is a paratopological group too \cite[Ex.1.9]{Ra2}. If
$(G,\tau)$ is a paratopological group then $\tau_r$ is the
strongest regular semigroup topology on the group $G$ which is
weaker than $\tau$; moreover, for any neighborhood base $\mathcal
B$ at the unit of the group $(G,\tau)$ the family $\mathcal
B_r=\{\inte\ol U:U\in\mathcal B\}$ is a base at the unit of the
group $(G,\tau_r)$ \cite[p.31--32]{Ra3}. The following proposition
is quite easy and probably is known.

\begin{proposition}\label{regularization} Let $(X,\tau)$ be a topological space. Then $(X,\tau)$ is
pseudocompact if and only if the regularization $(X,\tau_r)$ is
pseudocompact.
\end{proposition}

For the proof of Theorem~\ref{t3} we shall need a special
pseudocompact functionally Hausdorff semigroup topology on the
unit circle. We recall that a topological space $X$ is {\em
functionally Hausdorff\/} if continuous functions separate points
of $X$.

\begin{proposition}\label{pseudo} There is a functionally Hausdorff pseudocompact
first countable semigroup topology $\theta$ on the unit circle
$\IT$ which is not a  group topology.
\end{proposition}
\begin{proof}
Let $\IT$ be the unit circle and $\chi:\IT\to\IQ$ be a
(discontinuous) group homomorphism onto the groups of rational
numbers. Fix any element $x_0\in\IT$ with $\chi(x_0)=1$ and
observe that $S=\{1\}\cup\{x\in\IT:\chi(x)>0\}$ is a subsemigroup
of $\IT$. Let $\theta$ be the weakest semigroup topology on $\IT$
containing the standard compact topology $\tau$ and such that $S$
is open in $\theta$. It is easy to see that $\theta$ is
functionally Hausdorff and the sets $S\cap U$, where $1\in
U\in\tau$, form a neighborhood base of the topology $\theta$ at
the unit of $\IT$.

By Proposition~\ref{regularization}, to show that the group $(\IT,\theta)$ is
pseudocompact it suffices to verify that $\theta_r=\tau$. Since $\tau$ is a
regular semigroup topology on the group $\IT$ weaker than $\theta$, we get
$\theta_r\supset\tau$. To verify the inverse inclusion we first show that
${\overline U}^\tau={\overline U}^\theta$ for any $U\in\theta$. Since
$\tau\subset\theta$ it suffices to show that ${\overline
U}^\tau\subset{\overline U}^\theta$. Fix any point $x\in{\overline U}^\tau$ and
a neighborhood $V\in\tau$ of 1. We have to show that $x(V\cap S)\cap
U\ne\emptyset$. Pick up any point $y\in xV\cap U$. Since $U$ is open in the
topology $\theta$, we can find a neighborhood $W\in\tau$ of 1 such that
$y(W\cap S)\subset xV\cap U$. Find a number $N$ such that
$\chi(yx_0^N)>\chi(x)$ and thus $yx_0^n\in xS$ for all $n\ge N$ (we recall that
$x_0$ is an element of $\IT$ with $\chi(x_0)=1$). Moreover, since $x_0$ is
non-periodic in $\IT$, there exists a number $n\ge N$ such that $x_0^n\subset
W$. Then $yx_0^n\in (yS\cap yW)\cap xS\subset (xV\cap U)\cap xS=x(V\cap S)\cap
U$. Hence $x\in {\overline U}^\theta$ and
$\overline{U}^\theta=\overline{U}^\tau$.

Then  $$\inte_\theta{\ol U}^\theta=\IT \bs {\ol {\IT\bs{\ol
U}^\theta}}^\theta= \IT \bs {\ol {\IT\bs{\ol
U}^\theta}}^\tau\in\tau$$ which just yields $\theta_r\subset\tau$.
\end{proof}

Now we are able to present a {\em proof of Theorem~\ref{t3}}. The
``if'' part follows from the observation that for any Hausdorff
pseudocompact paratopological group $(G,\tau)$ its group reflexion
$G^\flat=(G,\tau_r)$ is a Hausdorff pseudocompact (and hence
totally bounded) topological group \cite{RR}.

To prove the ``only if'' part, fix a Bohr-separated abelian
paratopological group $(H,\tau)$ and let $\mathcal B_\tau$ be a
neighborhood base at the unit of the group $(H,\tau)$. It follows
that there is a group topology $\sigma'\subset\tau$ on $H$ such
that $(H,\sigma')$ is totally bounded. Let $(\hat H,\sigma)$ be
the Raikov completion of the group $(H,\sigma')$. It is clear that
$\hat H$ is a compact abelian group and $H$ is a normal dense
subgroup of $\hat H$. It follows that $\mathcal B_\tau$ is a
neighborhood base at the unit of some semigroup topology $\tau'$
on the group $\hat H$ with $\tau'|H=\tau$. Let $(\IT,\theta)$ be
the group from Proposition~\ref{pseudo}.

We shall define the topology on the product $G=\hat H\times \IT$
as follows. Let $\mathcal B_\tau$, $\mathcal B_\sigma$ and
$\mathcal B_\theta$ be the open neighborhood bases at the unit of
the groups $(H,\tau)$, $(\hat H,\sigma)$ and $(\IT,\theta)$
respectively. For arbitrary neighborhoods $U_\tau\in \mathcal
B_\tau$, $U_\sigma\in \mathcal B_\sigma$ and $U_\theta\in \mathcal
B_\theta$ with $U_\tau\subset U_\sigma$ let $[U_\tau,
U_\sigma,U_\theta]=U_\tau\times\{e_\IT\}\cup U_\sigma\times
(U_\theta\bs\{e_\IT\})$, where $e_H$ and $e_\IT$ are the units of
the groups $H$ and $\IT$ respectively. Denote by $\mathcal B$ the
family of all such $[U_\tau, U_\sigma,U_\theta]$. Repeating the
argument of Theorem~\ref{t1} check that the family $\mathcal B$ is
a base of some Hausdorff semigroup topology $\rho$ on $G$. By
$\pi$ denote the topology of the product $(\hat
H,\sigma)\times(\IT,\theta_r)$. By
Proposition~\ref{regularization} to show that the group $(G,\rho)$
is pseudocompact it suffice to verify that $\rho_r\subset\pi$.
For this we shall show that ${\ol U}^\rho\supset
U_\sigma\times{\ol U_\theta}^\theta$ for every $U=[U_\tau,
U_\sigma,U_\theta]\in \mathcal B$. Let $(x_{\hat H},x_\IT)\in
U_\sigma\times{\ol U_\theta}^\theta$ and $V=[V_\tau,
V_\sigma,V_\theta]\in \mathcal B$. It suffice to show that
$\big((x_{\hat
H},x_\IT)+V_\sigma\times(V_\theta\bs\{e_\IT\})\big)\cap
U_\sigma\times(U_\theta\bs\{e_\IT\})\not=\0$. This intersection is
nonempty if and only if the intersections $(x_{\hat
H}+V_\sigma)\cap U_\sigma$ and $(x_\IT+(V_\theta\bs\{e_\IT\}))\cap
(U_\theta\bs \{e_\IT\})$ are nonempty. The first intersection is
nonempty since $x_{\hat H}\in U_\sigma$ and the second is nonempty
since $x_\IT\in{\ol U_\theta}^\theta$ and the topology $\theta$ is
non-discrete.

\end{document}